\documentclass[10pt]{article}
\usepackage{latexsym}
\usepackage{amssymb}
\begin{document}
\title{\textbf {\small HOM-LIE TRIPLE SYSTEM AND HOM-BOL ALGEBRA STRUCTURES  ON  HOM-MALCEV AND RIGHT HOM-ALTERNATIVE 
ALGEBRAS}}
\author{{\sc Sylvain Attan and A. Nourou Issa}}
\date{}
\maketitle
\begin{abstract}
Every multiplicative Hom-Malcev algebra has a natural multiplicative Hom-Lie triple system structure. Moreover, there
is a natural Hom-Bol algebra structure on every multiplicative Hom-Malcev algebra and on every multiplicative right
(or left) Hom-alternative algebra. 
\end{abstract}
\footnote{2010 {\it Mathematics Subject Classification.} 17D10, 17D99
\par
{\it Key words and phrases.} Lie triple system, Malcev algebra, Hom-alternative algebra, Hom-Lie triple system, 
Hom-Malcev algebra, Hom-Bol algebra.} \\
\\
{\bf 1. Introduction} \\
\par
The study of Lie triple systems (Lts) on their own as algebraic objects started from Jacobson's work \cite{Jac} 
and developed further by, e.g., Lister \cite{Lister}, Yamaguti \cite{Yam} and other mathematicians. The 
interplay between Lts and the differential geometry of symmetric spaces is now folk (see, e.g., \cite{KN}, 
\cite{Loos2}). Lts constitute examples of ternary algebras. If $(g,[,])$ is a Lie algebra, then $(g, [,,])$ is a Lts, 
where $[x,y,z]:=[[x,y],z]$ (see \cite{Jac}, \cite{KN}, \cite{Loos2}). Another construction of Lts from binary algebras 
is the one from Malcev algebras found by Loos \cite{Loos1}.
\par
Malcev algebras were introduced by Mal'tsev \cite{Malc} in a study of commutator algebras of alternative 
algebras and also as a study of tangent algebras to local smooth Moufang loops. Mal'tsev used the name 
"Moufang-Lie algebras" for these nonassociative algebras while Sagle \cite{Sagle} introduced the term "Malcev 
algebras". Equivalent defining identities of Malcev algebras are pointed out in \cite{Sagle}.
\par
Alternative algebras, Malcev algebras and Lts (among other algebras) received a twisted generalization in the development of the theory 
of Hom-algebras during these latest years. The forerunner of the theory of Hom-algebras is the Hom-Lie algebra 
introduced by Hartwig, Larsson and Silvestrov in \cite{HLS} in order to describe the structure of some deformation of the Witt algebra and the Virasoro
 algebra. It is well-known that Lie algebras are related to associative algebras via the commutator bracket 
construction. In the search of a similar construction for Hom-Lie algebras, the notion of a Hom-associative 
algebra is introduced by Makhlouf and Silvestrov in \cite{MS1}, where it is proved that a Hom-associative algebra gives rise to a Hom-Lie 
algebra via the commutator bracket construction. Since then, various Hom-type structures are considered (see, e.g.,
\cite{AMS}, \cite{AI}, \cite{FGS}-\cite{Gohr},  \cite{Issa1}, \cite{Makh1}-\cite{MS2}, \cite{Yau3}-\cite{Yau6}). Roughly speaking, Hom-algebraic 
structures are corresponding ordinary algebraic structures whose defining identities are twisted by a linear 
self-map. A general method for constructing a Hom-type algebra from the ordinary type of algebra with a linear 
self-map is given by Yau in \cite{Yau2}.
\par
In \cite{AMS}, \cite{Yau4}, $n$-ary Hom-algebra structures generalizing n-ary algebras of Lie type or associative type 
were considered. In particular, generalizations of $n$-ary Nambu or Nambu-Lie algebras, called $n$-ary Hom-Nambu 
and Hom-Nambu-Lie algebras respectively, were introduced in \cite{AMS} while Hom-Jordan were defined in \cite{Makh2} 
and Hom-Lie triple systems 
(Hom-Lts) were introduced  in \cite{Yau4} (here, another definition of a Hom-Jordan algebra is given). It is shown 
\cite{Yau4} that Hom-Lts are ternary Hom-Nambu algebras with 
additional properties, that Hom-Lts arise also from Hom-Jordan triple systems or from other Hom-type algebras.
\par
Motivated by the relationships between some classes of binary algebras and some classes of binary-ternary algebras,
a study of Hom-type generalization of binary-ternary algebras is initiated in \cite{Issa1} with the definition of 
Hom-Akivis algebras. Further, Hom-Lie-Yamaguti algebras are considered in \cite{GI} and Hom-Bol algebras \cite{AI} 
are defined as a twisted generalization of Bol algebras which are introduced and studied in \cite{Mikh1}, \cite{SM1},
\cite{SM2} as infinitesimal structures tangent to smooth Bol loops (some aspects of the theory of Bol algebras are
discussed in \cite{Bouet}, \cite{HP} and \cite{PerIz}.
\par
In this paper, we will be concerned with right (or left) Hom-alternative algebras, Hom-Malcev algebras and Hom-Bol
algebras. We extend the Loos' construction of Lts from Malcev algebras (\cite{Loos1}, Satz 1) to the 
Hom-algebra setting (Section 3). Specifically, we prove (Theorem 3.2) that every multiplicative Hom-Malcev algebra 
is naturally a multiplicative Hom-Lts by a suitable definition of the ternary operation. As a tool in the proof 
of this fact, we point out a kind of compatibility relation between the original binary operation of a given 
Hom-Malcev algebra and the ternary operation mentioned above (Lemma 3.1). Moreover, we obtain that every multiplicative
Hom-Malcev algebra has a natural Hom-Bol algebra structure (Theorem 3.5). In \cite{Mikh2} Mikheev proved that every 
right alternative algebra has a natural (left) Bol algebra structure. In \cite{HP} Hentzel and Peresi proved that
not only a right alternative algebra but also a left alternative algebra have left Bol algebra structure. In 
Section 4 we prove that the Hom-analogue of these results hold. Specifically, every multiplicative right (or left)
Hom-alternative algebra is a Hom-Bol algebra (Theorem 4.4). In Section 2 we recall some basic definitions and facts
about Hom-algebras. We define the Hom-Jordan associator of a given Hom-algebra and point out that every Hom-algebra
is a Hom-triple system with respect to the Hom-Jordan associator. This observation is used in the proof of Theorem 4.4.   
\par
All vector spaces and algebras are meant over an algebraically closed ground field $\mathbb{K}$ of characteristic 0. \\
\\
{\bf 2. Some basics on Hom-algebras}\\
\par
We first recall some relevant definitions about binary and ternary Hom-algebras. In particular, we recall the 
notion of a Hom-Malcev algebra as well as some of its equivalent defining identities. Although various types of 
$n$-ary Hom-algebras are introduced and discussed in \cite{AMS}, \cite{Yau4}, for our purpose, we will consider 
ternary Hom-algebras (ternary Hom-Nambu algebras and Hom-Lts) and Hom-Bol algebras. For fundamentals on Hom-algebras, one may 
 refer, e.g., to \cite{AMS}, \cite{FGS}, \cite{HLS}, \cite{Makh1}, \cite{MS1}, \cite{Yau1}, \cite{Yau2}. Some aspects
of the theory of binary Hom-algebras are considered in \cite{Sheng}, while some classes of binary-ternary Hom-algebras
are defined and discussed in \cite{AI}, \cite{GI}, \cite{Issa1}.\\
\\
{\bf Definition 2.1.} (i) A {\it Hom-algebra} is a triple $(A, \ast, \alpha)$ in which $A$ is a $\mathbb{K}$-vector space, 
$\ast: A\times A\longrightarrow A$ a bilinear map (the binary operation) and  $\alpha: A\longrightarrow A$ a linear
map (the twisting map). The Hom-algebra $A$ is said to be {\it multiplicative} if $\alpha(x\ast y)=\alpha(x)\ast\alpha(y)$ 
for all $x,y \in A.$
\par 
(ii) The {\it Hom-Jacobian} in $(A,\ast, \alpha)$ is the trilinear map $J_{\alpha}: A \times A  \times A \longrightarrow A$  
defined as 
$J_{\alpha}(x,y,z):=\circlearrowleft_{x,y,z}(x\ast y)\ast\alpha(z),$  where $\circlearrowleft_{x,y,z}$ 
denotes the sum over cyclic permutation of $x,y,z.$
\par
(iii) The {\it Hom-associator} of a Hom-algebra $(A,\ast, \alpha)$ is the trilinear map $as: A^{\otimes 3}\longrightarrow A$
defined as $as(x,y,z) = (x \ast y) \ast \alpha (z) - \alpha (x) \ast (y \ast z)$. If $as(x,y,z) = 0$ for all $x,y,z \in A$,
then $(A, \ast, \alpha)$ is said to be {\it Hom-associative}.\\
\\
{\it Remark} 2.2. If $\alpha=id$ (the identity map), then a Hom-algebra $(A,\ast , \alpha)$ reduces to an ordinary 
algebra $(A,\ast)$, the Hom-Jacobian $J_{\alpha}$ is the ordinary Jacobian $J$, and the Hom-associator is the usual
associator for the algebra $(A,\ast).$\\
\par
As for ordinary algebras, to each Hom-algebra ${\cal A} := (A,\ast, \alpha)$ are attached two Hom-algebras: the
{\it commutator Hom-algebra} ${\cal A}^{-} := (A, [,], \alpha)$, where $[x,y] := x \ast y - y \ast x$ (the 
commutator of $x$ and $y$), and the {\it plus Hom-algebra} ${\cal A}^{+} := (A, \circ , \alpha)$, where $x \circ y :=
x \ast y + y \ast x$ (the {\it Jordan product}) for all $x,y \in A.$
\par
For our purpose, we make the following \\
\\
{\bf Definition 2.3.} The {\it Hom-Jordan associator} of a Hom-algebra ${\cal A} := (A,\ast, \alpha)$ is the trilinear
map $as^{J}: A^{\otimes 3}\longrightarrow A$ defined as  $as^{J}(x,y,z) = (x \circ y) \circ \alpha (z) - \alpha (x) \circ (y \circ z)$,
where "$\circ$" is the Jordan product on $A$.\\
\par
If $\alpha = id$, the Hom-Jordan associator reduces to the usual Jordan associator. \\
\\
{\bf Definition 2.4.} (i) A {\it Hom-Lie algebra} is  a Hom-algebra $(A, \ast, \alpha)$ such that the binary operation 
"$\ast$" is anticommutative and the {\it Hom-Jacobi identity} \\
\\
(2.1)  $J_{\alpha}(x,y,z)=0$ \\
\\
holds for all $x,y,z$ in $A$ (\cite{HLS}).
\par
(ii) {\it A Hom-Malcev algebra} (or {\it Hom-Maltsev algebra}) is a Hom-algebra $(A, \ast, \alpha)$ such that the 
binary operation "$\ast$" is anticommutative and that the {\it Hom-Malcev identity} \\
\\
(2.2)  $J_{\alpha}(\alpha(x),\alpha(y),x\ast z)=J_{\alpha}(x,y,z)\ast\alpha^2(x)$ \\
\\
holds for all $x,y,z$ in $A$ (\cite{Yau3}).
\par
(iii) A {\it Hom-Jordan algebra} is a Hom-algebra $(A, \ast, \alpha)$ such that $(A, \ast)$ is a commutative algebra
and the {\it Hom-Jordan identity}\\
\par
$as(x \ast x, \alpha (y), \alpha (x)) = 0$ \\
\\
is satisfied for all $x,y$ in $A$ (\cite{Yau3}).
\par
(iv) A Hom-algebra $(A, \ast, \alpha)$ is called a {\it right Hom-alternative algebra} if \\
\par
$as(x,y,y) = 0$ \\
\\
for all $x,y$ in $A$. A Hom-algebra $(A, \ast, \alpha)$ is called a {\it left Hom-alternative algebra} if \\
\par
$as(x,x,y) = 0$ \\
\\
for all $x,y$ in $A$. A Hom-algebra $(A, \ast, \alpha)$ is called a {\it Hom-alternative algebra} if it is both right
and left Hom-alternative (\cite{Makh2}). \\
\\
{\it Remark} 2.5. When $\alpha=id,$ the Hom-Jacobi identity (2.1) is the usual {\it Jacobi identity} $J(x,y,z)=0.$ 
Likewise, for $\alpha=id, $ the Hom-Malcev identity (2.2) reduces to the {\it Malcev identity} $J(x,y,x \ast z)=J(x,y,z) \ast x$. 
Therefore a Lie (resp. Malcev) algebra $(A,\ast)$ may be seen as a Hom-Lie (resp. Hom-Malcev) algebra with the 
identity map as the twisting map. Also Hom-Malcev algebras generalize  Hom-Lie algebras in the same way as Malcev 
algebras generalize Lie algebras. For $\alpha = id$
in the Hom-Jordan identity, we recover the usual Jordan identity. Observe that the definition of the Hom-Jordan 
identity in \cite{Yau3} is slightly different of the one formerly given in \cite{Makh2}.\\
\par
Hom-Malcev algebras are introduced in \cite{Yau3} in connection with a study of Hom-alternative algebras introduced 
in \cite{Makh2}. In fact it is proved (\cite{Yau3}, Theorem 3.8) that every Hom-alternative algebra is {\it Hom-Malcev 
admissible}, i.e. the commutator Hom-algebra of any Hom-alternative algebra is a Hom-Malcev algebra (this is the 
Hom-analogue of Mal'tsev's construction of Malcev algebras as commutator algebras of alternative algebras 
\cite{Malc}). This result is also mentioned in \cite{Issa1}, section 4, using an approach via Hom-Akivis algebras (this 
approach is close to the one of Mal'tsev in \cite{Malc}). Also, every Hom-alternative algebra is {\it Hom-Jordan
admissible}, i.e. its plus Hom-algebra is a Hom-Jordan algebra (\cite{Yau3}). Examples of Hom-alternative algebras 
and Hom-Jordan algebras could be found in \cite{Makh2} and \cite{Yau3}. An example of a right Hom-alternative 
algebra that is not left Hom-alternative is given in \cite{Yau6}.\\
\par
Equivalent to (2.2) defining identities of Hom-Malcev algebras are found in \cite{Yau3} where, in particular, it is 
shown that the identity \\
\\
(2.3) $J_{\alpha}(\alpha(x),\alpha(y),w\ast z)+J_{\alpha}(\alpha(w),\alpha(y),x\ast z)=$
\par
\hspace{3.5truecm}   $J_{\alpha}(x,y,z)\ast\alpha^2(w) +J_{\alpha}(w,y,z)\ast\alpha^2(x)$ \\
\\
is equivalent to (2.2) in any anticommutative Hom-algebra $(A, \ast, \alpha)$ (\cite{Yau3}, Proposition 2.7). In 
\cite{Issa2}, it is proved that in any anticommutative Hom-algebra $(A, \ast, \alpha)$, the Hom-Malcev identity (2.2) is 
equivalent to \\
\\
(2.4) $J_{\alpha}(\alpha(x),\alpha(y),u\ast v)=\alpha^2(u)\ast J_{\alpha}(x,y,v)+J_{\alpha}(x,y,u)\ast\alpha^2(v)$
\par
\hspace{3.5truecm} $-2J_{\alpha}(\alpha(u),\alpha(v),x\ast y)$.\\
\\
Therefore, apart from (2.2), the identities (2.3) and (2.4) may be taken as defining identities of a Hom-Malcev algebra.\\
\par
The Hom-algebras mentioned above are {\it binary} Hom-algebras. The first generalization of binary algebras was the
 ternary algebras introduced in \cite{Jac}. Ternary algebraic structures also appeared in various domains of 
theoretical and mathematical physics (see, e.g., \cite{Namb}). Likewise, 
binary Hom-algebras are generalized to $n$-ary Hom-algebra structures in \cite{AMS} (see also \cite{Yau4}).\\
\par
{\bf Definition 2.6} (\cite{AMS}). A {\it ternary Hom-Nambu algebra} is a triple $(A, [,,], \alpha)$ in which $A$ is 
a $\mathbb{K}$-vector space, $[,,]: A\times A\times A \longrightarrow A$ a trilinear map, and $\alpha=(\alpha_1, \alpha_2)$ 
a pair of linear maps (the twisting maps) such that the identity \\
\\
(2.5) $[\alpha_1(x),\alpha_2(y),[u,v,w]]=[[x,y,u],\alpha_1(v),\alpha_2(w)]+[\alpha_1(u),[x,y,v],\alpha_2(w)]$
\par
\hspace{4truecm}  $+[\alpha_1(u),\alpha_2(v),[x,y,w]]$ \\
\\
holds for all $u,v,w,x,y$ in $A.$ The identity (2.5) is called the {\it ternary Hom-Nambu identity}.\\
\\
{\it Remark} 2.7. When $\alpha_1=id=\alpha_2$ one recovers the usual ternary Nambu algebra. One may refer to \cite{Namb} 
for the origins of Nambu algebras. In \cite{AMS}, examples of $n$-ary Hom-Nambu algebras that are not Nambu algebras 
are provided.\\
\\ 
{\bf Definition 2.8} (\cite{Yau4}). A {\it Hom-Lie triple system} (Hom-Lts) is a ternary Hom-algebra $(A,[,,],\alpha=(\alpha_1,\alpha_2))$ 
such that \\
\\
(2.6) $[x,y,z]=-[y,x,z]$, \\
\\
(2.7) $\circlearrowleft_{x,y,z} [x,y,z]=0$, \\
\\
and the ternary Hom-Nambu identity (2.5) holds in $(A,[,,],\alpha=(\alpha_1,\alpha_2)).$ \\
\par
One notes that when the twisting maps $\alpha_1, \alpha_2$ are both equal to the identity map $id,$ then we 
recover the usual notion of a Lie triple system \cite{Lister}, \cite{Yam}. Examples of Hom-Lts could be found in
\cite{Yau4}. \\
\par
A particular situation , interesting for our setting, occurs when the twis-\\ting maps $\alpha_i$ are all equal, 
$\alpha_1=\alpha_2=\alpha$ and $\alpha([x,y,z])=[\alpha(x),\alpha(y),\alpha(z)]$ for all $x,y,z$ in $A.$ The 
Hom-Lie triple system $(A,[,,],\alpha)$ is then said to be {\it multiplicative} \cite{Yau4}. In case of 
multiplicativity, the ternary Hom-Nambu identity (2.5) then reads \\
\\
(2.8) $[\alpha(x),\alpha(y),[u,v,w]]=[[x,y,u],\alpha(v),\alpha(w)] +[\alpha(u),[x,y,v],\alpha(w)]$
\par
\hspace{3.5cm} $+[\alpha(u),\alpha(v),[x,y,w]]$. \\
\par
In \cite{GI} a (multiplicative) {\it Hom-triple system} is defined as a (multiplicative) ternary Hom-algebra 
$(A, [,,], \alpha)$ such that (2.6) and (2.7) are satisfied (thus a multiplicative Hom-Lts is seen as a Hom-triple 
system in which the identity (2.8) holds; observe that this definition of a Hom-triple system is different from 
the one formerly given in \cite{Yau4}, where a Hom-triple system is just the Hom-algebra $(A, [,,], \alpha)$).
With this vision of a Hom-triple system, it is shown (\cite{GI}) that every multiplicative non-Hom-associative algebra 
(i.e. not necessarily Hom-associative algebra) has a natural Hom-triple system structure if define $[x,y,z] := 
[[x,y], \alpha (z)] - as(x,y,z) + as(y,x,z)$. We note here that we get the same result if define another ternary
operation on a given Hom-algebra. Specifically, we have the following \\
\\
{\bf Proposition 2.9.} {\it Let ${\cal A} = (A, \ast , \alpha)$ be a multiplicative Hom-algebra. Define on 
$\cal A$ the ternary operation}\\
\\
(2.9) $(x,y,z) := as^{J}(y,z,x)$ \\
\\
{\it for all $x,y,z \in A$. Then $(A, (,,) , \alpha)$ is a multiplicative Hom-triple system}.\\
\\
{\it Proof.} A proof follows from the straightforward checking of the identities (2.6) and (2.7) for "$(,,)$" using
the commutativity of the Jordan product "$\circ$". \hfill $\square$ \\
\\
Since our results here depend on multiplicativity, in the rest of this paper we assume that all Hom-algebras
(binary or ternary) are multiplicative and while dealing with the binary operation "$\ast$" and where there is no 
danger of confusion, we will use juxtaposition in order to reduce the number of braces i.e., e.g., $xy\ast\alpha(z)$ 
means $(x\ast y)\ast \alpha(z).$
\par
Various results and constructions related to Hom-Lts are given in \cite{Yau4}. In particular, it is shown that every Lts $L$ 
can be twisted along any self-morphism of $L$ into a multiplicative Hom-Lts. For our purpose we just mention the 
following result. \\
\\
{\bf Proposition 2.10} (\cite{Yau4}). {\it Let $(A,\ast)$ be a Malcev algebra and $\alpha: A \longrightarrow A$ an 
algebra morphism. Then $A_{\alpha}:=(A,[,,]_{\alpha},\alpha)$ is a multiplicative Hom-Lts, where $[x,y,z]_{\alpha}=\alpha(2xy\ast z-yz\ast x-zx\ast y),$ 
for all $x,y,z$ in $A$}. \\
\par
One observes that the product $[x,y,z]=2xy\ast z-yz\ast x-zx\ast y$ is the one defined in \cite{Loos1} providing a 
Malcev algebra $(A,\ast)$ with a Lts structure. A  construction describing another view of Proposition 2.10 above, 
will be given in Section 3 (see Proposition 3.4) via Hom-Malcev algebras. For the time being, we point out the 
following slight generalization of the result above, producing a sequence of multiplicative Hom-Lts from a given 
Malcev algebra.\\
\\
{\bf Proposition 2.11.} {\it Let $(A,\ast)$ be a Malcev algebra and $\alpha: A \longrightarrow A$ an algebra 
morphism. Let $\alpha^{0} = id$ and, for any integer $n \geq 1$, $\alpha^{n} = \alpha \circ  \alpha^{n-1}$. If 
define on $A$ a trilinear operation $[ , , ]_{\alpha^{n}}$ by
\par
$ [x,y,z]_{\alpha^{n}} = \alpha^{n} (2xy \ast z - yz \ast x - zx \ast y)$ \\
for all $x,y,z$ in $A$, then $(A,[,,]_{\alpha^{n}},\alpha^{n})$ is a multiplicative Hom-Lts}. \\
\\
{\it Proof.} Let $[x,y,z] = 2xy\ast z-yz\ast x-zx\ast y$ and then $ [x,y,z]_{\alpha^{n}} = \alpha^{n} ([x,y,z])$.
We shall use the fact that $(A, [ , , ])$ is a Lts (\cite{Loos1}). The identities (2.6) and (2.7) for $[x,y,z]_{\alpha^{n}}$ 
are quite obvious. Next,
\par
$[\alpha^{n} (x), \alpha^{n} (y) , [u,v,w]_{\alpha^{n}}]_{\alpha^{n}} = [\alpha^{n} (x), \alpha^{n} (y) , \alpha^{n} ([u,v,w])]_{\alpha^{n}} $
\par
$= \alpha^{2n} ([x,y,[u,v,w]])$
\par
$ = \alpha^{2n} ([[x,y,u],v,w]) + \alpha^{2n} ([u,[x,y,v],w]) + \alpha^{2n} ([u,v,[x,y,w]])$
\par
$= [\alpha^{n} ([x,y,u]), \alpha^{n} (v), \alpha^{n} (w)]_{\alpha^{n}} + [\alpha^{n} (u), \alpha^{n} ([x,y,v]), \alpha^{n} (w) ]_{\alpha^{n}}$
\par
$+ [\alpha^{n} (u), \alpha^{n} (v) , \alpha^{n} ([x,y,w]) ]_{\alpha^{n}}$
\par
$= [[x,y,u]_{\alpha^{n}}, \alpha^{n} (v), \alpha^{n} (w)]_{\alpha^{n}} + [\alpha^{n} (u), [x,y,v]_{\alpha^{n}},\alpha^{n} (w) ]_{\alpha^{n}}$
\par
$ + [\alpha^{n} (u), \alpha^{n} (v) ,[x,y,w]_{\alpha^{n}}]_{\alpha^{n}} $ \\
and so (2.8) holds for $[ , , ]_{\alpha^{n}}$. Thus $(A,[,,]_{\alpha^{n}},\alpha^{n})$ is a multiplicative Hom-Lts. \hfill $\square$ \\
\par
In \cite{AI} we defined a Hom-Bol algebra as a twisted generalization of a (left) Bol algebra. For the introduction
and original studies
of Bol algebras, we refer to \cite{Mikh1}, \cite{SM1}, \cite{SM2}. Bol algebras are further considered in, e.g., 
\cite{HP}, \cite{PerIz}. \\
\\
{\bf Definition 2.12} (\cite{AI}). A {\it Hom-Bol algebra} is a quadruple $(A, [,], (,,), \alpha)$ in which $A$ is
a vector space, "$[,]$" a binary operation, "$(,,)$" a ternary operation on $A$, and $\alpha : A \rightarrow A$ a linear 
map such that\\
(HB1) $\alpha ([x,y]) = [\alpha (x), \alpha (y)]$, \\
(HB2) $\alpha ((x,y,z)) = (\alpha (x),\alpha (y),\alpha (z))$, \\
(HB3) $[x,y] = - [y,x]$, \\
(HB4) $(x,y,z) = - (y,x,z)$, \\
(HB5) ${\circlearrowleft}_{x,y,z}(x,y,z) = 0$,\\
(HB6) $(\alpha (x),\alpha (y),[u,v]) = [(x,y,u), {\alpha}^{2} (v)] + [{\alpha}^{2} (u), (x,y,v)]$ 
\par
\hspace{3.5truecm} $+ (\alpha (u),\alpha (v),[x,y]) - [[\alpha (u),\alpha (v)], [\alpha (x),\alpha (y)]]$,\\
(HB7) $({\alpha}^{2} (x), {\alpha}^{2} (y), (u,v,w)) = ((x,y,u), {\alpha}^{2} (v), {\alpha}^{2} (w)) + 
{\alpha}^{2} (u), (x,y,v), {\alpha}^{2} (w))$
\par 
\hspace{4.5truecm} $+ ({\alpha}^{2} (u), {\alpha}^{2} (v), (x,y,w))$\\
for all $u,v,w,x,y,z \in A$. \\
\par
The identities (HB1) and (HB2) mean the multiplicativity of $(A, [,], (, ,), \alpha)$. It is built into our definition
for convenience.
\par
One observes that for $\alpha = id$ the identities (HB3)-(HB7) reduce to the defining identities of a (left) Bol
algebra \cite{Mikh1} (see also \cite{HP}, \cite{PerIz}). If $[x,y] = 0$ for all $x,y \in A$, then $(A, [,], (, ,), \alpha)$
becomes a (multiplicative) Hom-Lts $(A, (, ,), {\alpha}^{2})$.
\par
Construction results and some examples of Hom-Bol algebras are given in \cite{AI}. In particular, Hom-Bol algebras
can be constructed from Malcev algebras. The Hom-analogues of the construction of Bol algebras from Malcev algebras
\cite{Mikh1} or from right alternative algebras \cite{Mikh2} (see also \cite{HP}) are considered in this paper. \\
\\
{\bf 3. Hom-Lts and Hom-Bol algebras from Hom-Malcev algebras}\\
\par
In this section, we prove that every multiplicative Hom-Malcev algebra has a natural multiplicative Hom-Lts structure
(Theorem 3.2) and, moreover, a natural Hom-Bol algebra structure (Theorem 3.5). Theorem 3.2  could be seen as the 
Hom-analogue of the Loos' result (\cite{Loos1}, Satz 1). Besides the identities (2.3) and (2.4), Lemma 3.1 below is a 
tool in the proof of this result. Theorem 3.5 could be seen as the Hom-analogue of a construction by Mikheev \cite{Mikh1}
of Bol algebras from Malcev algebras. Proposition 3.4 is another view of a result in \cite{Yau4} (see Proposition 2.10
above).
\par
In his work \cite{Loos1}, Loos considered in a Malcev algebra $(A,\ast),$ the following ternary operation: \\
\\
(3.1) $\{x,y,z\}= 2xy\ast z-yz\ast x-zx\ast y$. \\
\\
Then $(A,\{,,\})$ turns out to be a Lts. This result, in the Hom-algebra setting, looks as in Theorem 3.2 below. 
Similarly as in the Loos construction, our investigations are based on the following ternary operation in a 
Hom-Malcev algebra $(A,\ast,\alpha)$: \\
\\
(3.2) $\{x,y,z\}_{\alpha}=2xy\ast\alpha(z)-yz\ast\alpha(x)-zx\ast\alpha(y).$ \\
\\
From (3.2) it clearly follows that $\{,,\}_{\alpha}$ can also be written as \\
\\
(3.3) $\{x,y,z\}_{\alpha}=-J_{\alpha}(x,y,z)+3xy \ast \alpha(z)$. \\
\\
One observes that when $\alpha=id,$ we recover the product (3.1). This agrees with the reduction of the Hom-Malcev 
algebra $(A,\ast,\alpha)$ to the  Malcev algebra $(A,\ast)$. First, we prove the following\\
\\
{\bf Lemma 3.1.} {\it Let $(A,\ast,\alpha)$ be a Hom-Malcev algebra. If define on $(A,\ast,\alpha)$ a ternary 
operation $"\{,,\}_{\alpha}"$ by} (3.2), {\it then} \\
\\
(3.4) $\{\alpha(x),\alpha(y) , u\ast v\}_{\alpha} = \alpha^2(u)\ast\{x,y,v\}_{\alpha} +\{x,y,u\}_{\alpha} \ast \alpha^2(v)$
\par
\hspace{4cm} $-J_{\alpha}(\alpha(u),\alpha(v), x\ast y)$ \\
\\
{\it for all $u, v,x,y$ in $A$}.\\
\\
{\it Proof.} Let us write (2.4) as
\par
$-J_{\alpha}(\alpha(x),\alpha(y),u\ast v)= -J_{\alpha}(x,y,u)\ast\alpha^2(v) +\alpha^2(u)\ast(-J_{\alpha}(x,y,v))$
\par
\hspace{3.5cm} $+3J_{\alpha}(\alpha(u),\alpha(v),x\ast y)-J_{\alpha}(\alpha(u),\alpha(v),x\ast y)$ \\
i.e.
\par
$-J_{\alpha}(\alpha(x),\alpha(y),u\ast v) = -J_{\alpha}(x,y,u)\ast\alpha^2(v) +\alpha^2(u)\ast(-J_{\alpha}(x,y,v))$
\par
\hspace{3.5cm} $+3\alpha(u)\alpha(v)\ast\alpha(x\ast y)+3(\alpha(v)\ast xy)\ast\alpha^2(u)$
\par
\hspace{3.5cm}                                       $+3(xy\ast\alpha(u))\ast\alpha^2(v)-J_{\alpha}(\alpha(u),\alpha(v),x\ast y)$. \\
Therefore, by multiplicativity, we have
\par
$-J_{\alpha}(\alpha(x),\alpha(y),u\ast v)+3\alpha(x)\alpha(y)\ast\alpha(u\ast v)$
\par
$=(-J_{\alpha}(x,y,u)+3xy\ast\alpha(u))\ast\alpha^2(v)+\alpha^2(u)\ast(-J_{\alpha}(x,y,v)+3xy\ast\alpha(v))$
\par
$-J_{\alpha}(\alpha(u),\alpha(v),x\ast y)$ \\
and so, we get (3.4) by (3.3). \hfill $\square$\\
\\
We now prove \\
\\
{\bf Theorem 3.2.} {\it Let $(A,\ast,\alpha)$ be a multiplicative Hom-Malcev algebra. If define on $(A,\ast,\alpha)$ 
a ternary operation $"\{,,\}_{\alpha}"$ by} (3.2), {\it then $(A,\{,,\}_{\alpha},\alpha^2)$ is a multiplicative Hom-Lts.}\\
\\
{\it Proof.} We must prove the validity of (2.6), (2.7), (2.8) for the operation (3.2) in the Hom-Malcev algebra $(A,\ast,\alpha)$.
\par
First observe that the multiplicativity of $(A,\ast,\alpha)$ implies that $\alpha^2 (\{x,y,z\}_{\alpha}) \\ = \{\alpha^2(x),\alpha^2(y),\alpha^2(z) \}_{\alpha}$, 
with $x,y,z$ in $A$.
\par
From the skew-symmetry of $"\ast"$ and $J_{\alpha}(x,y,z),$ it clearly follows from (3.3) that  $\{x,y,z\}_{\alpha} = -\{y,x,z\}_{\alpha}$ 
which is (2.6) for $"\{,,\}_{\alpha}".$
\par
Next, using (3.3) and the skew-symmetry of $J_{\alpha}(x,y,z)$ where applicable, we compute
\par
$\{x,y,z\}_{\alpha}+\{y,z,x\}_{\alpha}+\{z,x,y\}_{\alpha}$
\par
$= -J_{\alpha}(x,y,z)+3xy\ast\alpha(z)-J_{\alpha}(y,z,x)+3yz\ast\alpha(x)-J_{\alpha}(z,x,y)+3zx\ast\alpha(y)$
\par
$=-3J_{\alpha}(x,y,z)+3J_{\alpha}(x,y,z) =0$ \\
and thus $\circlearrowleft_{x,y,z}\{x,y,z\}_{\alpha}=0,$ so we get (2.7) for "$\{,,\}_{\alpha}$".
\par
Consider now $\{\alpha^2(x),\alpha^2(y),\{u,v,w\}_{\alpha}\}_{\alpha}$ in $(A,\ast,\alpha)$. Then
\begin{eqnarray}
&&\{\alpha^2(x),\alpha^2(y),\{u,v,w\}_{\alpha}\}_{\alpha}\nonumber\\
&=&\{\alpha^2(x),\alpha^2(y),2uv\ast\alpha(w)-vw\ast\alpha(u)-wu\ast\alpha(v)\}_{\alpha} \mbox{ (by (3.2))} \nonumber\\
&=&\{\alpha^2(x),\alpha^2(y),2uv\ast\alpha(w)\}_{\alpha}-\{\alpha^2(x),\alpha^2(y),vw\ast\alpha(u)\}_{\alpha} \nonumber\\
 &&- \{\alpha^2(x),\alpha^2(y),wu\ast\alpha(v)\}_{\alpha}\nonumber\\
&=&\{\alpha(x),\alpha(y),2u\ast v\}_{\alpha}\ast\alpha^3(w)+\alpha^2(2u\ast v)\ast\{\alpha(x),\alpha(y),\alpha(w)\}_{\alpha}\nonumber\\
&&-J_{\alpha}(\alpha(2u\ast v),\alpha^2(w),\alpha(x\ast y))-\{\alpha(x),\alpha(y),v\ast w\}_{\alpha}\ast\alpha^3(u)\nonumber\\
&&-\alpha^2(v\ast w)\ast\{\alpha(x),\alpha(y),\alpha(u)\}_{\alpha}+J_{\alpha}(\alpha(v\ast w),\alpha^2(u),\alpha(x\ast y))\nonumber\\
&&-\{\alpha(x),\alpha(y),w\ast u\}_{\alpha}\ast\alpha^3(v)-\alpha^2(w\ast u)\ast\{\alpha(x),\alpha(y),\alpha(v)\}_{\alpha}\nonumber\\
&&+J_{\alpha}(\alpha(w\ast u),\alpha^2(v),\alpha(x\ast y)) \mbox{  (by (3.4))}\nonumber\\
&=&(2\{x,y,u\}_{\alpha}\ast\alpha^2(v)+2\alpha^2(u)\ast\{x,y,v\}_{\alpha}\nonumber\\
&&-2J_{\alpha}(\alpha(u),\alpha(v),x\ast y))\ast\alpha^3(w)+2\alpha^2(u\ast v)\ast\{\alpha(x),\alpha(y),\alpha(w)\}_{\alpha}\nonumber\\
&&-J_{\alpha}(\alpha(2u\ast v),\alpha^2(w),\alpha(x\ast y))-(\{x,y,v\}_{\alpha}\ast\alpha^2(w)\nonumber\\
&&+\alpha^2(v)\ast\{x,y,w\}_{\alpha}-J_{\alpha}(\alpha(v),\alpha(w),x\ast y))\ast\alpha^3(u)\nonumber\\
&&-\alpha^2(v\ast w)\ast\{\alpha(x),\alpha(y),\alpha(u)\}_{\alpha}+J_{\alpha}(\alpha(v\ast w),\alpha^2(u),\alpha(x\ast y))\nonumber\\
&&-(\{x,y,w\}_{\alpha}\ast\alpha^2(u)+\alpha^2(w)\ast\{x,y,u\}_{\alpha} \nonumber\\
&&- J_{\alpha}(\alpha(w),\alpha(u),x\ast y))\ast\alpha^3(v)-\alpha^2(w\ast u)\ast\{\alpha(x),\alpha(y),\alpha(v)\}_{\alpha}\nonumber\\
&&+J_{\alpha}(\alpha(w\ast u),\alpha^2(v),\alpha(x\ast y))\mbox{   ( again by (3.4)) } \nonumber\\
&=&2\{x,y,u\}_{\alpha}\alpha^2(v)\ast\alpha^3(w)+2\alpha^2(u)\{x,y,v\}_{\alpha}\ast\alpha^3(w)\nonumber\\
&&-2J_{\alpha}(\alpha(u),\alpha(v),x\ast y)\ast\alpha^3(w)+2\alpha^2(u\ast v)\ast\{\alpha(x),\alpha(y),\alpha(w)\}_{\alpha}\nonumber\\
&&-J_{\alpha}(\alpha(2u\ast v),\alpha^2(w),\alpha(x\ast y))-\{x,y,v\}_{\alpha}\alpha^2(w)\ast\alpha^3(u)\nonumber\\
&&-\alpha^2(v)\{x,y,w\}_{\alpha}\ast\alpha^3(u)+J_{\alpha}(\alpha(v),\alpha(w),x\ast y)\ast\alpha^3(u)\nonumber\\
&&-\alpha^2(v\ast w)\ast\{\alpha(x),\alpha(y),\alpha(u)\}_{\alpha}+J_{\alpha}(\alpha(v\ast w),\alpha^2(u),\alpha(x\ast y))\nonumber\\
&&-\{x,y,w\}_{\alpha}\alpha^2(u)\ast\alpha^3(v)-\alpha^2(w)\{x,y,u\}_{\alpha}\ast\alpha^3(v)\nonumber\\
&&+J_{\alpha}(\alpha(w),\alpha(u),x\ast y)\ast\alpha^3(v)-\alpha^2(w\ast u)\ast\{\alpha(x),\alpha(y),\alpha(v)\}_{\alpha}\nonumber\\
&&+J_{\alpha}(\alpha(w\ast u),\alpha^2(v),\alpha(x\ast y))\nonumber\\
&=&2\{x,y,u\}_{\alpha}\alpha^2(v)\ast\alpha^3(w)-\alpha^2(v\ast w)\ast\alpha(\{x,y,u\}_{\alpha}\nonumber\\
&&-\alpha^2(w)\{x,y,u\}_{\alpha}\ast\alpha^3(v)+2\alpha^2(u)\{x,y,v\}_{\alpha}\ast\alpha^3(w)\nonumber\\
&&-\{x,y,v\}_{\alpha}\alpha^2(w)\ast\alpha^3(u)-\alpha^2(w\ast u)\ast\alpha(\{x,y,v\}_{\alpha})\nonumber\\
&&+2\alpha^2(u\ast v)\ast\alpha(\{x,y,w\}_{\alpha})-\alpha^2(v)\{x,y,w\}_{\alpha}\ast\alpha^3(u)\nonumber\\
&&-\{x,y,w\}_{\alpha}\alpha^2(u)\ast\alpha^3(v)-2J_{\alpha}(\alpha(u),\alpha(v),x\ast y)\ast\alpha^3(w)\nonumber\\
&&-J_{\alpha}(\alpha(2u\ast v),\alpha^2(w),\alpha(x\ast y))+J_{\alpha}(\alpha(v),\alpha(w),x\ast y)\ast\alpha^3(u)\nonumber\\
&&+J_{\alpha}(\alpha(v\ast w),\alpha^2(u),\alpha(x\ast y))+J_{\alpha}(\alpha(w),\alpha(u),x\ast y)\ast\alpha^3(v)\nonumber\\
&&+J_{\alpha}(\alpha(w\ast u),\alpha^2(v),\alpha(x\ast y))\nonumber\\
 &&\mbox{  (rearranging terms)}\nonumber\\
&=&\{\{x,y,u\}_{\alpha},\alpha^2(v),\alpha^2(w)\}_{\alpha}+\{\alpha^2(u),\{x,y,v\}_{\alpha},\alpha^2(w)\}_{\alpha}\nonumber\\
&&+\{\alpha^2(u),\alpha^2(v),\{x,y,w\}_{\alpha}\}_{\alpha} \nonumber \\
&&+\mbox{ \bf [ } -2(J_{\alpha}(\alpha(u),\alpha(v),x\ast y) \ast\alpha^3(w)
+J_{\alpha}(\alpha(u\ast v),\alpha^2(w),\alpha(x\ast y))) \nonumber\\
&&+J_{\alpha}(\alpha(v),\alpha(w),x\ast y)\ast\alpha^3(u)
+J_{\alpha}(\alpha(v\ast w),\alpha^2(u),\alpha(x\ast y)) \nonumber\\
&&+J_{\alpha}(\alpha(w),\alpha(u),x\ast y)\ast\alpha^3(v)
+J_{\alpha}(\alpha(w\ast u),\alpha^2(v),\alpha(x\ast y)) \mbox{\bf ] }. \nonumber
\end{eqnarray}
In this latest expression, denote by $N(u,v,w,x,y)$ the expression in "{\bf [}...{\bf ]}". To conclude, we proceed 
to show that $N(u,v,w,x,y)=0.$
\par
Observe first that, by (2.3), we have
\begin{eqnarray}
&&J_{\alpha}(\alpha(u),x\ast y,\alpha(w))\ast\alpha^2(\alpha(v))+J_{\alpha}(\alpha(v),x\ast y,\alpha(w))\ast\alpha^2(\alpha(u))\nonumber\\
&=& J_{\alpha}(\alpha^2(u),\alpha(x\ast y),\alpha(v)\ast\alpha(w))+J_{\alpha}(\alpha^2(v),\alpha(x\ast y),\alpha(u)\ast\alpha(w))\nonumber
\end{eqnarray}
i.e.,
\begin{eqnarray}
&&J_{\alpha}(\alpha(w),\alpha(u),x\ast y)\ast\alpha^3(v)+J_{\alpha}(\alpha(w\ast u),\alpha^2(v),\alpha(x\ast y))\nonumber\\
&=& J_{\alpha}(\alpha (v\ast w),\alpha^2(u),\alpha(x\ast y))+J_{\alpha}(\alpha(v),\alpha(w),x\ast y)\ast\alpha^3(u)\nonumber
\end{eqnarray}
With this observation, the expression $N(u,v,w,x,y)$ is transformed as follows:
\begin{eqnarray}
&&N(u,v,w,x,y)\nonumber\\
&=&2[-J_{\alpha}(\alpha(u),\alpha(v),x\ast y)\ast\alpha^3(w)-J_{\alpha}(\alpha(u\ast v),\alpha^2(w),\alpha(x\ast y))]\nonumber\\
&&+2[J_{\alpha}(\alpha(v\ast w),\alpha^2(u),\alpha(x\ast y))+J_{\alpha}(\alpha(v),\alpha(w),x\ast y)\ast\alpha^3(u)]\nonumber\\
&=&2[-J_{\alpha}(\alpha^2(w),\alpha(x\ast y),\alpha(u)\ast\alpha(v))-J_{\alpha}(\alpha^2(u),\alpha(x\ast y),\alpha(w)\ast\alpha(v))\nonumber\\
&&-J_{\alpha}(\alpha(u),\alpha(v),x\ast y)\ast\alpha^3(w)+J_{\alpha}(\alpha(v),\alpha(w),x\ast y)\ast\alpha^3(u)]\nonumber\\
&=&2[-J_{\alpha}(\alpha(w),x\ast y,\alpha(v))\ast\alpha^3(u))-J_{\alpha}(\alpha(u),x\ast y, \alpha(v))\ast\alpha^3(w)\nonumber\\
&&-J_{\alpha}(\alpha(u),\alpha(v),x\ast y)\ast\alpha^3(w)+J_{\alpha}(\alpha(v),\alpha(w),x\ast y)\ast\alpha^3(u)]\nonumber\\
&&\mbox{ (applying (2.3) to} -J_{\alpha}(\alpha^2(w),\alpha(x\ast y), \alpha(u)\ast\alpha(v)) \nonumber \\ &&-J_{\alpha}(\alpha^2(u),\alpha(x\ast y),\alpha(w)\ast\alpha(v)) )\nonumber \\
&=&0  \mbox{ ( by the skew-symmetry of } J_{\alpha}(x,y,z) ). \nonumber
\end{eqnarray}
Therefore, we obtain that (2.8) holds for "$\{,,\}_{\alpha}$" and thus $(A,\{,,\}_{\alpha},\alpha^2)$ is a Hom-Lts.
 \hfill $\square$\\
\\
{\it Remark} 3.3.  In the proof of his result, Loos (\cite{Loos1}, Satz 1) used essentially the fact that the left 
translations $L(x)$ in a Malcev algebra $(A,\ast)$ are derivations of the ternary operation "$\{,,\}$" defined by 
(3.1). Unfortunately, for Hom-Malcev algebras such a tool is still 
not available at hand.\\
\par
From \cite{Yau3} (Theorem 2.12) we know that any Malcev algebra $A$ can be twisted into a Hom-Malcev algebra along 
any linear self-map of $A.$ Consistent with this result, we recall the following method for constructing Hom-Lts 
which, in fact, is a result in \cite{Yau4} (see also Proposition 2.10 and Proposition 2.11 above) but using a Hom-Malcev 
algebra construction in our proof (as a consequence of Theorem 3.2).\\
\\
{\bf Proposition 3.4.} {\it Let $(A,\ast)$ be a Malcev algebra and $\alpha$ any self-morphism of $(A,\ast)$. If 
define on $(A,\ast)$ a ternary operation "$\{,,\}_{\alpha}$" by}
\begin{eqnarray}
\{x,y,z\}_{\alpha}&=&\alpha^2(2xy\ast z-yz\ast x-zx\ast y),\nonumber
\end{eqnarray}
{\it then $(A,\{,,\}_{\alpha},\alpha^2)$ is a multiplicative Hom-Lts}.\\
\\
{\it Proof.} One knows (\cite{Yau3}, Theorem 2.12) that from $(A,\ast)$ and any self-morphism $\alpha$ of $(A,\ast)$, 
we get a (multiplicative) Hom-Malcev algebra $(A,\widetilde{\ast},\alpha)$, where $x\widetilde{\ast}y=\alpha(x\ast y)$ for all $x,y$ in $A$.
Next, if define on $(A,\widetilde{\ast},\alpha)$ a ternary operation
\begin{eqnarray}
\{x,y,z\}_{\alpha}&:=&2(x\widetilde{\ast}y)\widetilde{\ast}\alpha(z)-(y\widetilde{\ast}z)\widetilde{\ast}\alpha(x)
-(z\widetilde{\ast}x)\widetilde{\ast}\alpha(y),\nonumber
\end{eqnarray}
then by Theorem 3.2, $(A,\{,,\}_{\alpha},\alpha^2)$ is a Hom-Lts and "$\{,,\}_{\alpha}$" expresses through "$\ast$" as
\begin{eqnarray}
\{x,y,z\}_{\alpha} &=& 2\alpha(\alpha(x\ast y)\ast\alpha(z))-\alpha(\alpha(y\ast z)\ast\alpha(x)) -\alpha(\alpha(z\ast x)\ast\alpha(y))\nonumber\\
&=&2\alpha^2(xy\ast z)-\alpha^2(yz\ast x)-\alpha^2(zx\ast y)\nonumber\\
&=&\alpha^2(2xy\ast z-yz\ast x-zx\ast y).  \nonumber
\end{eqnarray}
\hfill $\square$ \\
\par
Observe that, although constructed in a quite different way, the operation "$\{,,\}_{\alpha}$" in Proposition 3.4 
above coincide with "$[ , , ]_{{\alpha}^{n}}$" in Proposition 2.11 for $n=2$.
\par
Combining Lemma 3.1 and Theorem 3.2, we get the following result. \\
\\
{\bf Theorem 3.5.} {\it Let  $(A,\ast , \alpha)$ be a multiplicative Hom-Malcev algebra. If define on  $(A,\ast , \alpha)$ 
a ternary operation $(,,)_{\alpha}$ by} \\
\\
(3.5) $(x,y,z)_{\alpha} := \frac{1}{3} \{x,y,z \}_{\alpha}$,\\
\\
{\it where "$\{,,\}_{\alpha}$" is defined by (3.3), then $(A,\ast , (,,)_{\alpha}, \alpha)$ is a Hom-Bol algebra}.\\
\\
{\it Proof}. The definition (3.5) and Theorem 3.2 imply that $(A, (,,)_{\alpha}, {\alpha}^{2})$ is a multiplicative
Hom-Lts i.e. (HB2), (HB4), (HB5) and (HB7) hold for $(A,\ast , (,,)_{\alpha}, \alpha)$. Now, (HB1) and (HB3) are
respectively the multiplicativity and skew-symmetry of "$\ast$". Next, we are done if we prove (HB6) for 
$(A,\ast , (,,)_{\alpha}, \alpha)$.
\par
From (3.3) and multiplicativity we have \\
$-J_{\alpha}(\alpha(u), \alpha(v), x \ast y) = \{ \alpha(u), \alpha(v), x \ast y \}_{\alpha} -3( \alpha(u)\alpha(v)) \ast 
(\alpha(x)\alpha(y))$ \\
and then (3.4) takes the form \\
$ \{ \alpha(x), \alpha(y), u \ast v \}_{\alpha} = \{x,y,u \}_{\alpha} \ast \alpha^2(v) + \alpha^2(u) \ast \{x,y,v \}_{\alpha}$
\par 
\hspace{3truecm} $ + \{ \alpha(u), \alpha(v), x \ast y \}_{\alpha} -3( \alpha(u)\alpha(v)) \ast (\alpha(x)\alpha(y))$. \\
\\
Multiplying by $\frac{1}{3}$ each member of this latter equality and using (3.5), we get \\
$ ( \alpha(x), \alpha(y), u \ast v )_{\alpha} = (x,y,u )_{\alpha} \ast \alpha^2(v) + \alpha^2(u) \ast (x,y,v )_{\alpha}$
\par 
\hspace{3truecm} $ + ( \alpha(u), \alpha(v), x \ast y )_{\alpha} -( \alpha(u)\alpha(v)) \ast (\alpha(x)\alpha(y))$ \\
which is (HB6) for $(A,\ast , (,,)_{\alpha}, \alpha)$. So $(A,\ast , (,,)_{\alpha}, \alpha)$ is a Hom-Bol algebra. \hfill $\square$\\
\par
Since any Hom-alternative algebra is Hom-Malcev admissible (\cite{Yau3}, Theorem 3.8), from Theorem 3.5 we have the
following \\
\\
{\bf Corollary 3.6.} {\it Let $(A, \ast, \alpha)$ be a multiplicative Hom-alternative algebra. Then $(A, [,], (,,)_{\alpha}, \alpha)$ is
a Hom-Bol algebra, where $(x,y,z)_{\alpha} := -\frac{1}{3} J_{\alpha}(x,y,z)$ \\ $+ xy \ast \alpha(z)$, for 
all $x,y \in A$}. \hfill $\square$\\
\par
The aim of Section 4 is a generalization of Corollary 3.6 to multiplicative right (or left) Hom-alternative algebras. \\
\par
Various constructions of Hom-Lts are offered in \cite{Yau4} starting from either Hom-associative algebras, Hom-Lie 
algebras, Hom-Jordan triple systems, ternary totally Hom-associative algebras, Malcev algebras or alternative algebras. 
In practice, it is easier to construct Hom-Lts or Hom-Bol algebras from well-known (binary) algebras such as, e.g., alternative 
algebras or Malcev algebras. From this point of view, our construction results (Theorem 3.2, Proposition 3.4 and
Theorem 3.5) have rather a theoretical feature (the extension to Hom-algebra setting of the Loos' result 
\cite{Loos1} and a result by Mikheev \cite{Mikh1}) than a practical method for constructing Hom-Lts or Hom-Bol algebras.
However, it could be of some interest to get a Hom-Lts or a Hom-Bol algebra from a given Hom-Malcev algebra without
resorting to the corresponding Malcev algebra.\\
\\
{\bf 4. Hom-Lts and Hom-Bol algebras from right (or left) Hom- alternative algebras}\\
\par
In this section we prove that every multiplicative right (or left) Hom-alternative algebra has a natural Hom-Bol
algebra structure (and, subsequently, a natural Hom-Lts structure). This is the Hom-analogue of a result by Mikheev
\cite{Mikh2} and by Hentzel and Peresi \cite{HP}.
\par
First we recall some few basic properties of right Hom-alternative algebras that could be found in \cite{Makh2},
\cite{Yau6}.
\par
The linearized form of the right Hom-alternative identity $as(x,y,y) =0$ is given by the following result. \\
\\
{\bf Lemma 4.1} (\cite{Makh2}). {\it If $(A, \ast, \alpha)$ is a Hom-algebra, then the following statements are 
equivalent. 
\par
(i) $(A, \ast, \alpha)$ is right Hom-alternative.
\par
(ii) $(A, \ast, \alpha)$ satisfies} \\
\\
(4.1) $as(x,y,z) = - as(x,z,y)$ \\
\\
{\it for all $x,y,z \in A$.
\par
(iii) $(A, \ast, \alpha)$ satisfies} \\
\\
(4.2) $\alpha (x) \ast (yz + zy) = xy \ast \alpha (z) + xz \ast \alpha (y)$\\
\\
{\it for all $x,y,z \in A$}. \\
\par
Observe that if $(A, \ast, \alpha)$ is a right Hom-alternative algebra, then $(A, {\ast}^{op}, \alpha)$ is a left
Hom-alternative algebra, where $x {\ast}^{op} y := y \ast x$. So the mirrors of (4.1) and (4.2) hold for $(A, {\ast}^{op}, \alpha)$: \\
\\
(4.3) $as(x,y,z) = - as(y,x,z)$ \\
\\
and\\
\\
(4.4) $((x {\ast}^{op} y) + (y {\ast}^{op} x)) {\ast}^{op} \alpha (z) = \alpha (x) {\ast}^{op} (y {\ast}^{op} z) 
+ \alpha (y) {\ast}^{op} (x {\ast}^{op} z)$.\\
\par
Now we have the following\\
\\
{\bf Lemma 4.2.} {\it In any multiplicative right Hom-alternative algebra $(A, \ast, \alpha)$, the identity}\\
\\
(4.5) $as([u,v], \alpha (x), \alpha (y)) =[as(u,x,y),{\alpha}^{2}(v)] + [{\alpha}^{2}(u), as(v,x,y)] $
\par
\hspace{3.5truecm} $+ as(\alpha (v), \alpha (u), [x,y]) -  as(\alpha (u), \alpha (v), [x,y])$ \\
\\
{\it holds for all} $x,y,z \in A$.\\
\\
{\it Proof}. The identity
\par
$ as(uv, \alpha (x), \alpha (y)) = as(u,x,y){\alpha}^{2}(v) + {\alpha}^{2}(u)as(v,x,y) - as(\alpha (u), \alpha (v), [x,y])$\\
holds in any right Hom-alternative algebra (see \cite{Yau6}, Theorem 7.1 (7.1.1c)). Next, in this identity, switching
$u$ and $v$, we have
\par
$ as(vu, \alpha (x), \alpha (y)) = as(v,x,y){\alpha}^{2}(u) + {\alpha}^{2}(v)as(u,x,y) - as(\alpha (v), \alpha (u), [x,y])$.\\
Then, subtracting memberwise this latter equality from the one above and using the linearity of $as$, we get (4.5). \hfill $\square$\\
\par
Note that in case when $(A, \ast, \alpha)$ is a left Hom-alternative algebra, the identity (4.5) reads as \\
\\
(4.6)  $as(\alpha (x), \alpha (y), [u,v]) =[as(x,y,u),{\alpha}^{2}(v)] + [{\alpha}^{2}(u), as(x,y,v)] $
\par
\hspace{3.5truecm} $+ as( [x,y],\alpha (v), \alpha (u)) -  as([x,y],\alpha (u), \alpha (v))$. \\
\par
In any multiplicative right (or left) Hom-alternative algebra $(A, \ast, \alpha)$ we consider the ternary operation
defined by (2.9), i.e. \\
\par
$(x,y,z) := as^{J}(y,z,x)$, \\
\\
where $ as^{J}$ is the Hom-Jordan associator defined in Section 2. Observe that for $\alpha = id$ the ternary operation
"$(,,)$" is precisely the one defined in \cite{HP} (see also \cite{Mikh2}, Remark 2) and that makes any right (or left)
alternative algebra into a left Bol algebra. In \cite{HP}, Hentzel and Peresi used the approach of Mikheev \cite{Mikh2}
who formerly proved that the commutator algebra of any right alternative algebra has a left Bol algebra structure.\\
\\
{\bf Proposition 4.3.} {\it (i) If $(A, \ast, \alpha)$ is a multiplicative right Hom-alternative algebra, then}\\
\\
(4.7) $(x,y,z) = [[x,y],\alpha (z)] - 2as(z,x,y)$\\
\\
{\it for all} $x,y,z \in A$.
\par
{\it (ii) If $(A, \ast, \alpha)$ is a multiplicative left Hom-alternative algebra, then}\\
\\
(4.8) $(x,y,z) = [[x,y],\alpha (z)] - 2as(x,y,z)$\\
\\
{\it for all} $x,y,z \in A$. \\
\\
{\it Proof}. (i) From (2.9) we have
\begin{eqnarray}
(x,y,z) &=& (y \circ z) \circ \alpha (x) - \alpha (y) \circ (z \circ x) \nonumber \\
&=& ((y \ast z) + (z \ast y)) \ast \alpha (x) + [\alpha (x) \ast  ((y \ast z) + (z \ast y))]\nonumber \\
&-& [\alpha (y) \ast  ((z \ast x) + (x \ast z))] - ((z \ast x) + (x \ast z)) \ast \alpha (y)\nonumber \\
&=& ((y \ast z) + (z \ast y)) \ast \alpha (x) + [(x \ast y) \ast \alpha (z) + (x \ast z) \ast \alpha (y)]\nonumber \\
&-& [(y \ast z) \ast \alpha (x) + (y \ast x) \ast \alpha (z)] \nonumber \\
&-& ((z \ast x) + (x \ast z)) \ast \alpha (y) \;\; \mbox{(by (4.2))}\nonumber \\
&=&  (z \ast y) \ast \alpha (x) + (x \ast y) \ast \alpha (z) - (y \ast x) \ast \alpha (z)
- (z \ast x) \ast \alpha (y) \nonumber \\
&=& (z \ast y) \ast \alpha (x) - (z \ast x) \ast \alpha (y) + [x,y] \ast \alpha (z) \nonumber \\
&=& (z \ast y) \ast \alpha (x) - (z \ast x) \ast \alpha (y) + [[x,y], \alpha (z)] + \alpha (z) \ast [x,y] \nonumber \\
&=& [[x,y], \alpha (z)] + (z \ast y) \ast \alpha (x) - \alpha (z) \ast (y \ast x) - (z \ast x) \ast \alpha (y)\nonumber \\
&+& \alpha (z) \ast (x \ast y) \nonumber \\
&=& [[x,y], \alpha (z)] + as(z,y,x) - as(z,x,y) \nonumber \\
&=& [[x,y], \alpha (z)] - 2as(z,x,y) \;\; \mbox{(by (4.1))} \nonumber 
\end{eqnarray}
and so we get (4.7). \\
(ii) Proceeding as above but using (4.4) and then (4.3), one gets (4.8). \hfill $\square$ \\
\par
We are now in position to prove the main result of this section. \\
\\
{\bf Theorem 4.4.} {\it Let $(A, \ast, \alpha)$ be a multiplicative right (resp. left) Hom-alternative algebra. If
define on $A$ a ternary operation "$(,,)$" by (4.7) (resp. (4.8)), then $(A, (,,), {\alpha}^{2})$ is a Hom-Lts and 
$(A, [,], (,,), \alpha )$ is a Hom-Bol algebra}. \\
\\
{\it Proof}. We prove the theorem for a multiplicative right Hom-alternative algebra $(A, \ast, \alpha)$ (the proof
of the left case is the mirror of the right one).
\par
The identities (HB1) and (HB2) follow from the multiplicativity of $(A, \ast, \alpha)$. The identities (HB3) and
(HB4) are obvious from the definition of "$[,]$" and "$(,,)$". The identity (HB5) follows from Proposition 2.9.
\par
In \cite{Yau5} Yau showed that if, on a multiplicative Hom-Jordan algebra  $(A, \circ, \alpha)$, define a ternary
operation by \\
\par
$[x,y,z] := 2(\alpha (x) \circ (y \circ z) - \alpha (y) \circ (x \circ z)$, \\
\\
then $(A, [,,], {\alpha}^{2})$ is a multiplicative Hom-Lts (see \cite{Yau5}, Corollary 4.1). Now, observe that
$[x,y,z] = 2as^{J}(y,z,x)$, i.e. $[x,y,z] = 2(x,y,z)$. Therefore, since every multiplicative right Hom-alternative
algebra is Hom-Jordan admissible (see \cite{Yau6}, Theorem 4.3), we conclude that $(A, (,,), {\alpha}^{2})$ is a 
multiplicative Hom-Lts and so the identity (HB7) holds for $(A, [,], (,,), \alpha)$.
\par
Next, $(A, [,], (,,), \alpha )$ is a Hom-Bol algebra if we prove that (HB6) additionally holds.
\par
Write (4.7) as\\
\\
(4.9) $-2as(z,x,y) = (x,y,z) - [[x,y], \alpha (z)]$. \\
\\
Multiplying each member of (4.5) by $-2$ and next using (4.9), we get
\par
$(\alpha (x), \alpha (y), [u,v]) - [[\alpha (x), \alpha (y)], \alpha ([u,v])]$
\par
$= [(x,y,u)-[[x,y], \alpha (u)], {\alpha}^{2}(v)] + [{\alpha}^{2}(u), (x,y,v)-[[x,y], \alpha (v)]] $
\par
$+ (\alpha (u), [x,y], \alpha (v)) - [[\alpha (u),[x,y]], {\alpha}^{2}(v)]$
\par
$- (\alpha (v), [x,y], \alpha (u)) + [[\alpha (v),[x,y]], {\alpha}^{2}(u)]$\\
i.e. \\
\\
(4.10) $ (\alpha (x), \alpha (y), [u,v]) = [(x,y,u), {\alpha}^{2}(v)] + [{\alpha}^{2}(u), (x,y,v)]$
\par
\hspace{0.5truecm} $- ([x,y], \alpha (u), \alpha (v)) + ([x,y], \alpha (v), \alpha (u)) + \alpha ( [[x,y],[u,v]])$. \\
\\
Observe that $ -([x,y], \alpha (u), \alpha (v)) + ([x,y], \alpha (v), \alpha (u))$ \\
$= (\alpha (u), [x,y], \alpha (v)) + ([x,y], \alpha (v), \alpha (u))$ \\
$= - (\alpha (v), \alpha (u), [x,y])$ (since $ \circlearrowleft_{a,b,c} (a,b,c) =0$ by (HB5))\\
$= (\alpha (u), \alpha (v), [x,y])$.\\
Therefore, (4.10) now reads \\
$ (\alpha (x), \alpha (y), [u,v]) = [(x,y,u), {\alpha}^{2}(v)] + [{\alpha}^{2}(u), (x,y,v)]$
\par
\hspace{2.5truecm}$+(\alpha (u), \alpha (v), [x,y]) - \alpha ( [[u,v],[x,y]])$\\
and so (HB6) holds for $(A, [,], (,,), \alpha)$. Thus we conclude that $(A, [,], (,,), \alpha)$ is a Hom-Bol algebra.
One gets the same result in case when $(A, \ast, \alpha)$ is a multiplicative left Hom-alternative algebra and essentially 
using (4.8) and (4.6). This finishes the proof. \hfill $\square$

{\sc Sylvain Attan} \\
Institut de Math\'ematiques et de Sciences Physiques,  \\ Universit\'e d'Abomey-Calavi, 01 BP 613-Oganla Porto-Novo, 
B\'enin \\
{\it E-mail address}: sylvain.attan@imsp-uac.org \\ and \\
{\sc A. Nourou Issa}\\
D\'epartement de Math\'ematiques, Universit\'e d'Abomey-Calavi, 01 BP 4521 Cotonou, 
B\'enin \\
{\it E-mail address}: woraniss@yahoo.fr
\end{document}